\documentclass[11pt,a4paper]{article}
\usepackage[T1]{fontenc}
\usepackage[utf8]{inputenc}
\usepackage[english]{babel}
\usepackage{amsmath,amssymb,amsthm,mathtools}
\usepackage{newtxtext,newtxmath}
\usepackage[margin=24mm,top=20mm,bottom=22mm,headheight=22pt,headsep=16pt,footskip=24pt]{geometry}
\usepackage{microtype}
\usepackage{graphicx}
\usepackage{caption}
\usepackage{placeins}
\usepackage{xcolor}
\definecolor{ink}{HTML}{242930}
\definecolor{sectionblue}{HTML}{315E77}
\definecolor{questionback}{HTML}{F4F8FC}
\definecolor{questionframe}{HTML}{377EB8}
\definecolor{answerback}{HTML}{F4FAF5}
\definecolor{answerframe}{HTML}{45945B}
\definecolor{aiback}{HTML}{FFF9E9}
\definecolor{aiframe}{HTML}{D99A20}
\definecolor{aititle}{HTML}{F2CD72}
\usepackage{titlesec}
\titleformat{\section}{\large\bfseries\color{sectionblue}}{\thesection}{0.7em}{}
\titleformat{\subsection}{\normalsize\bfseries}{\thesubsection}{0.7em}{}
\titleformat{\subsubsection}{\normalsize\bfseries}{}{0pt}{}
\titlespacing*{\section}{0pt}{2.0ex plus .5ex minus .3ex}{0.8ex}
\titlespacing*{\subsection}{0pt}{1.8ex plus .4ex minus .2ex}{0.6ex}
\titlespacing*{\subsubsection}{0pt}{1.5ex plus .3ex minus .2ex}{0.4ex}
\usepackage{fancyhdr}
\renewcommand{\headrulewidth}{0.3pt}
\renewcommand{\headrule}{\hbox to\headwidth{\color{sectionblue!45}\leaders\hrule height\headrulewidth\hfill}}
\fancypagestyle{plain}{\fancyhf{}\fancyfoot[C]{\thepage}\renewcommand{\headrulewidth}{0pt}}
\usepackage[most]{tcolorbox}
\tcbset{enhanced,breakable,sharp corners,boxrule=0.45pt,
  left=7pt,right=7pt,top=6pt,bottom=6pt,before skip=9pt,after skip=9pt}
\usepackage{enumitem}
\setlist{topsep=4pt,itemsep=2pt,parsep=0pt,leftmargin=1.6em}
\usepackage{needspace}
\usepackage[unicode,colorlinks=true,linkcolor=sectionblue,citecolor=sectionblue,urlcolor=sectionblue]{hyperref}
\hypersetup{pdftitle={Separately Continuous Functions on Squares and Mixed Products},
  pdfsubject={Baire class one on squares and mixed products},
  pdfauthor={}}
\newtheorem{theorem}{Theorem}[section]
\newtheorem{lemma}[theorem]{Lemma}
\newtheorem{proposition}[theorem]{Proposition}
\newtheorem{corollary}[theorem]{Corollary}
\theoremstyle{definition}

\theoremstyle{remark}

\numberwithin{equation}{section}
\newcommand{\R}{\mathbb R}
\newcommand{\N}{\mathbb N}
\newcommand{\bb}{\mathfrak b}
\newcommand{\SC}{\operatorname{SC}}
\newcommand{\supp}{\operatorname{supp}}
\newcommand{\fin}{[\omega_1]^{<\omega}}

\makeatletter
\renewcommand{\l@section}{\@dottedtocline{1}{0em}{1.4em}}
\renewcommand{\l@subsection}{\@dottedtocline{2}{1.4em}{2.0em}}
\makeatother
\begin{document}
\color{ink}
\thispagestyle{plain}
\begin{center}
{\fontsize{18}{21}\selectfont\bfseries On Separately Continuous Functions of the First\\[0.10em]
Baire Class on Squares and Mixed Products\par}
\vspace{0.45em}
{\normalsize\itshape A Conditional Negative Answer to Mykhaylyuk's Problem~42\par}

\vspace{0.8em}
{\normalsize Alexander V. Osipov\footnote{Email: \href{mailto:OAB@list.ru}{\texttt{OAB@list.ru}} .}\par}
\vspace{0.35em}
{\small\itshape Krasovskii Institute of Mathematics and Mechanics, Yekaterinburg, Russia\par}

\end{center}
\vspace{0.4em}
\begin{list}{}{\setlength{\leftmargin}{14pt}\setlength{\rightmargin}{14pt}}
\item[]\small\textbf{\color{sectionblue}Abstract.}
We consider Mykhaylyuk's problem asking whether the Baire-one property
for every separately continuous real-valued function on each of the
squares $X^2$ and $Y^2$ implies the same property on $X\times Y$.
Assuming $\mathfrak b>\omega_1$, the article obtains a negative answer
using the compact cube $K=[-1,1]^{\omega_1}$ and the space
$E=(c_{00}(\omega_1),\sigma(c_{00}(\omega_1),\ell_\infty(\omega_1)))$.
It presents a finite-coordinate characterization of Baire-one functions
on $E$ under this hypothesis. 
\end{list}
\vspace{0.4em}
\noindent{\small\textbf{\color{sectionblue}Keywords.} Separate continuity; Baire class one; weak topology; coordinate approximation; countable chain condition; bounding number.\par}
\vspace{0.45em}
\noindent{\small\textbf{\color{sectionblue}2020 Mathematics Subject Classification.} 54C08, 54C30, 03E17.\par}
\vspace{0.65em}
\tableofcontents
\section*{Introduction}
\addcontentsline{toc}{section}{Introduction}

We consider Problem~42 of V.~Mykhaylyuk (January~7, 2022),
recorded in the Lviv Scottish Book; see~\cite{ScottishBook}.
\begin{tcolorbox}[colback=questionback,colframe=questionframe]
\noindent\textbf{Question (Mykhaylyuk, Problem 42).}
Let $X$ and $Y$ be completely regular spaces such that
\[
\SC(X^2,\R)\subseteq B_1(X^2),\qquad
\SC(Y^2,\R)\subseteq B_1(Y^2).
\]
Must $\SC(X\times Y,\R)\subseteq B_1(X\times Y)$ hold?
\end{tcolorbox}
Here $\SC(S\times T,\R)$ denotes the class of separately continuous
real-valued functions, and $B_1(Z)$ denotes the class of pointwise
limits of sequences of continuous real-valued functions on $Z$.

For $u,v\in\N^\N$, write $u\leq^*v$ if $u(n)\leq v(n)$ for all
sufficiently large $n$. The bounding number $\bb$ is the least
cardinality of an unbounded family in $(\N^\N,\leq^*)$.
Throughout the main construction we assume
\begin{equation}\label{assumption}
 \omega_1<\bb.
\end{equation}
Consequently, every family of at most $\omega_1$ members of
$\N^\N$ has a common upper bound with respect to $\leq^*$.

Put
\[
 V=c_{00}(\omega_1),\qquad
 E=(V,\sigma(V,\ell_\infty(\omega_1))).
\]
The elements of $V$ have finite support. The weak topology above
is generated by the linear functionals
\[
 L_b(a)=\sum_{\alpha<\omega_1}b_\alpha a_\alpha,
 \qquad b\in\ell_\infty(\omega_1).
\]
It is a Hausdorff locally convex topology; in particular, $E$ is
completely regular. For a finite set $A\subseteq\omega_1$, let
\[
 E_A=\{a\in E:\supp a\subseteq A\}.
\]
The topology induced on $E_A$ is its usual Euclidean topology.
Write $\pi_Aa$ for the vector obtained from $a$ by replacing all
coordinates outside $A$ by zero, and write $e_\alpha$ for the
$\alpha$th coordinate unit vector.

By an algebraic polynomial of finite degree on $V$ we mean a
finite sum of diagonal restrictions of algebraic multilinear
forms and a constant. Equivalently, such a polynomial is a
formal sum of coordinate monomials with a common finite bound
on their degrees and arbitrary real coefficients.
Only finitely many monomials are nonzero at any given vector of
$V$. No continuity of these polynomials is assumed.

\section{A preliminary reduction via a topological sum}

Let $P(Z)$ denote the assertion that every separately continuous
real-valued function on $Z\times Z$ belongs to $B_1(Z^2)$.

\begin{proposition}\label{sumprop}
For arbitrary topological spaces $X$ and $Y$, the following conditions
are equivalent:
\begin{enumerate}[label=(\roman*)]
\item $P(X\sqcup Y)$;
\item $P(X)$, $P(Y)$, and
$\SC(X\times Y,\R)\subseteq B_1(X\times Y)$.
\end{enumerate}
\end{proposition}
\begin{proof}
The square of the disjoint topological sum $X\sqcup Y$
is the disjoint union of four clopen subspaces
$X^2$, $X\times Y$, $Y\times X$, and $Y^2$.
Assume (i). A separately continuous function on any of these
four products extends by zero to a separately continuous
function on $(X\sqcup Y)^2$; restricting its continuous
approximations establishes the required assertions in (ii).
The case $Y\times X$ follows by interchanging the factors.

Conversely, let $F$ be separately continuous on
$(X\sqcup Y)^2$. Its restriction to each of the four clopen
products has a sequence of continuous pointwise approximations
by (ii). At each stage, glue these four continuous functions
along the clopen partition. This gives a continuous function
on the whole square, and the resulting sequence converges
pointwise to $F$.
\end{proof}

Thus the original question is equivalent to asking whether
spaces with property $P$ are closed under finite topological sums.
This reformulation alone does not imply an affirmative answer.

\section{Approximation of multilinear forms}

\begin{lemma}\label{multilinear}
Under \eqref{assumption}, every algebraic $d$-linear form
$B:V^d\to\R$, where $d\geq1$, is the limit of a sequence of
jointly continuous $d$-linear forms on $E^d$, uniformly on every
product of compact subsets of finite-dimensional coordinate
subspaces.
\end{lemma}

\begin{proof}
Write
\[
 B(a^1,\ldots,a^d)=
 \sum_{\alpha_1,\ldots,\alpha_d<\omega_1}
 c_{\alpha_1,\ldots,\alpha_d}
 \prod_{j=1}^d a^j_{\alpha_j}.
\]
For fixed arguments this is a finite sum.

For each $\beta<\omega_1$, choose a bijective enumeration of
$\beta+1$ whose domain is a subset of $\N$.
Given a tuple $(\alpha_1,\ldots,\alpha_d)$, let $\beta$ be its
maximum and let $p$ be the first position at which this maximum
occurs. Fixing $p$ and the enumeration indices of the other
$\alpha_j$ partitions $\omega_1^d$ into at most countably many
sets, each of the form
\[
 \{(\gamma_1(\beta),\ldots,\gamma_d(\beta)):\beta\in D\},
 \qquad \gamma_p(\beta)=\beta,\quad
 \gamma_j(\beta)\leq\beta.
\]
The requirement that $p$ be the first position of the maximum
is included in the definition of $D$, so these sets are indeed
pairwise disjoint.

The corresponding part of the form is
\begin{equation}\label{graph}
 G(a^1,\ldots,a^d)=
 \sum_{\beta\in D}w_\beta a^p_\beta
       \prod_{j\ne p}a^j_{\gamma_j(\beta)}.
\end{equation}
Choose distinct numbers $t_\alpha\in(0,1)$, $\alpha<\omega_1$,
and put
\[
 k_m(s,t)=\cos^{2m}(\pi(s-t)).
\]
Then
\begin{equation}\label{kernel-limit}
 k_m(t_\alpha,t_\gamma)\longrightarrow
 \mathbf1_{\{\alpha=\gamma\}}.
\end{equation}
Indeed, for $\alpha\ne\gamma$ we have
$0<|t_\alpha-t_\gamma|<1$, and therefore
$|\cos(\pi(t_\alpha-t_\gamma))|<1$.
The cosine subtraction formula and the binomial theorem give
\[
 k_m(s,t)=\sum_{\nu=1}^{N_m}u_{m,\nu}(s)v_{m,\nu}(t),
\]
where all factors are bounded on $(0,1)$.

Let $w_\beta^{(m)}=\max(-m,\min(m,w_\beta))$, and define
\begin{equation}\label{graph-approx}
 \begin{split}
 G_m(a^1,\ldots,a^d)
 =\sum_{\beta\in D}w_\beta^{(m)}a^p_\beta
 \prod_{j\ne p}\left(
    \sum_{\alpha<\omega_1}a^j_\alpha
       k_m(t_\alpha,t_{\gamma_j(\beta)})
                    \right).
 \end{split}
\end{equation}
For fixed arguments all sums are finite.
Expanding the finite representations of the kernels expresses
$G_m$ as a finite sum of products of functionals of the form
\[
 \sum_\alpha a^j_\alpha u_{m,\nu_j}(t_\alpha)
 \quad(j\ne p)
\]
and a functional of the form
\[
 \sum_{\beta\in D}a^p_\beta w_\beta^{(m)}
       \prod_{j\ne p}v_{m,\nu_j}(t_{\gamma_j(\beta)}).
\]
All coefficient families are bounded. Hence these linear
functionals are weakly continuous, and $G_m$ is jointly
continuous on $E^d$.

For fixed arguments, the index $\beta$ in \eqref{graph-approx}
ranges over the finite set $D\cap\supp a^p$.
Thus \eqref{kernel-limit} and $w_\beta^{(m)}\to w_\beta$
imply that $G_m\to G$ pointwise.

Enumerate the parts as $G_1,G_2,\ldots$, appending zero forms
if there are only finitely many parts. Put
\[
 S_N=\sum_{j=1}^N G_j,\qquad
 S_{N,m}=\sum_{j=1}^N G_{j,m}.
\]
For every $\boldsymbol\alpha=(\alpha_1,\ldots,\alpha_d)$,
the coefficient $S_N(e_{\alpha_1},\ldots,e_{\alpha_d})$
equals $c_{\boldsymbol\alpha}$ for all sufficiently large $N$.
For each $N$, choose $r_{\boldsymbol\alpha}(N)\in\N$ such that
\[
 |S_{N,m}(e_{\alpha_1},\ldots,e_{\alpha_d})
       -S_N(e_{\alpha_1},\ldots,e_{\alpha_d})|<1/N
 \quad\text{whenever }m\geq r_{\boldsymbol\alpha}(N).
\]
Since $|\omega_1^d|=\omega_1<\bb$, there exists $m\in\N^\N$
eventually dominating every $r_{\boldsymbol\alpha}$.
It follows that the coefficients of $S_{N,m(N)}$ converge
to the corresponding coefficients of $B$.

Only finitely many coefficients occur on
$E_{A_1}\times\cdots\times E_{A_d}$ for finite sets $A_j$.
Their convergence implies uniform convergence on every product
of compact subsets of these subspaces.
\end{proof}

\begin{corollary}\label{polynomials}
Under \eqref{assumption}, every algebraic polynomial
$P:V\to\R$ of finite degree is the limit of a sequence of
continuous polynomials on $E$, uniformly on every compact subset
of every $E_A$, $A\in\fin$.
\end{corollary}
\begin{proof}
Express $P$ as a finite sum of diagonal restrictions of
multilinear forms and a constant, and apply
Lemma~\ref{multilinear} to each form.
\end{proof}

\section{A finite-dimensional criterion for Baire class one}

\begin{theorem}\label{criterion}
Under \eqref{assumption}, the following conditions are equivalent
for an arbitrary function $f:E\to\R$:
\begin{enumerate}
 \item $f\in B_1(E)$;
 \item $f|_{E_A}\in B_1(E_A)$ for every $A\in\fin$.
\end{enumerate}
\end{theorem}

\begin{proof}
The implication $(1)\Rightarrow(2)$ follows by restricting
continuous approximations. We prove the converse.

For $A\in\fin$ and $x\in E_A$, put
\[
 f_A(x)=\sum_{B\subseteq A}(-1)^{|A\setminus B|}f(\pi_Bx).
\]
In particular, $f_\varnothing=f(0)$. Inclusion--exclusion
inversion yields
\begin{equation}\label{mobius}
 f(x)=\sum_{A\subseteq\supp x}f_A(\pi_Ax).
\end{equation}
If $x_\alpha=0$ for some $\alpha\in A$, the terms defining
$f_A$ cancel in pairs, and hence $f_A(x)=0$.

By assumption, $f_A\in B_1(E_A)$.
Every Baire class one function on $\R^k$ is a pointwise limit
of polynomials: take continuous approximations $g_n$, and
approximate $g_n$ by a polynomial within $1/n$ on $[-n,n]^k$,
using the Weierstrass polynomial approximation theorem on
finite-dimensional cubes.

Choose polynomials converging pointwise to $f_A$ and apply
the operator
\[
 \Delta_Ap(x)=
 \sum_{B\subseteq A}(-1)^{|A\setminus B|}p(\pi_Bx).
\]
Since $f_A$ vanishes whenever any coordinate indexed by $A$
is zero, $\Delta_Af_A=f_A$.
We obtain polynomials $p_{A,n}\to f_A$ such that every monomial
of $p_{A,n}$ contains all the variables indexed by $A$.
For $A=\varnothing$, take $p_{\varnothing,n}=f(0)$.

For this argument, assign degree zero to the zero polynomial.
The family of sequences $n\mapsto\deg p_{A,n}$ has cardinality
at most $\omega_1$. By \eqref{assumption}, there is $D\in\N^\N$
such that, for every fixed $A$,
\[
 \deg p_{A,n}\leq D(n)
\]
for all sufficiently large $n$. Define
\[
 q_{A,n}=
 \begin{cases}
 p_{A,n},&\deg p_{A,n}\leq D(n),\\
 0,&\deg p_{A,n}>D(n),
 \end{cases}
 \qquad
 P_n(x)=\sum_{A\in\fin}q_{A,n}(\pi_Ax).
\]
For fixed $x$, a summand can be nonzero only if
$A\subseteq\supp x$. Moreover, all monomials have degree at
most $D(n)$. Thus $P_n$ is an algebraic polynomial of finite
degree on $V$. For each fixed $A$, we have $q_{A,n}\to f_A$.
It follows from \eqref{mobius} that
\begin{equation}\label{Pn-limit}
 P_n(x)\longrightarrow f(x)\qquad(x\in E).
\end{equation}

By Corollary~\ref{polynomials}, choose continuous polynomials
$P_{n,m}$ converging, as $m\to\infty$, to $P_n$ uniformly on
every compact subset of each finite-dimensional coordinate
subspace. For $A\in\fin$ and $r\in\N$, put
\[
 C_{A,r}=\{x\in E_A: |x_\alpha|\leq r\ (\alpha\in A)\}.
\]
Choose $s_{A,r}(n)$ so that
\[
 \sup_{x\in C_{A,r}}|P_{n,m}(x)-P_n(x)|<1/n
 \quad\text{whenever }m\geq s_{A,r}(n).
\]
The family of sequences $s_{A,r}$ has cardinality at most
$\omega_1$. Choose a common eventual upper bound $m(n)$.
Every $x\in E$ belongs to some $C_{A,r}$, so
\[
 |P_{n,m(n)}(x)-P_n(x)|<1/n
\]
for all sufficiently large $n$.
Together with \eqref{Pn-limit}, this gives
$P_{n,m(n)}(x)\to f(x)$, proving that $f\in B_1(E)$.
\end{proof}

\section{Separate continuity and a metrizable factor}

\begin{lemma}\label{metric}
If $M$ is metrizable, $T$ is a topological space, and
$F:T\times M\to\R$ is separately continuous, then
$F\in B_1(T\times M)$.
\end{lemma}
\begin{proof}
Fix a compatible metric on $M$. For each $n$, choose a locally
finite partition of unity $(\varphi_{n,i})_{i\in I_n}$
subordinate to an open cover $(U_{n,i})_{i\in I_n}$ by sets of
diameter less than $1/n$, and choose $z_{n,i}\in U_{n,i}$.
Such partitions exist by paracompactness of metrizable spaces.
Put
\[
 F_n(t,z)=\sum_{i\in I_n}\varphi_{n,i}(z)F(t,z_{n,i}).
\]
The sum is locally finite, so $F_n$ is continuous.
If $\varphi_{n,i}(z)>0$, then $z,z_{n,i}\in U_{n,i}$.
For fixed $t$, continuity of $F(t,\cdot)$ at $z$ therefore
implies $F_n(t,z)\to F(t,z)$.
\end{proof}

\begin{corollary}\label{Esquare}
Under \eqref{assumption},
$\SC(E^2,\R)\subseteq B_1(E^2)$.
\end{corollary}
\begin{proof}
The space $E\times E$ is linearly homeomorphic to $E$.
Indeed, the product of the two weak topologies corresponds
to the dual pair $(c_{00}(I),\ell_\infty(I))$, where $I$ is the
disjoint union of two copies of $\omega_1$ and hence has
cardinality $\omega_1$.

If $F:E\times E\to\R$ is separately continuous, then
Lemma~\ref{metric} gives
$F|_{E_A\times E_B}\in B_1(E_A\times E_B)$ for finite $A,B$.
Every finite-dimensional coordinate subspace of $E\times E$
has this form. Apply Theorem~\ref{criterion} after the
identification above.
\end{proof}

\begin{lemma}\label{compact}
If $K$ is a separable compact Hausdorff space and $L$ is a
compact Hausdorff space, then
$\SC(K\times L,\R)\subseteq B_1(K\times L)$.
\end{lemma}
\begin{proof}
Let $D\subseteq K$ be countable and dense, and let
$F:K\times L\to\R$ be separately continuous. Define
\[
 q:L\to\R^D,\qquad q(y)=(F(d,y))_{d\in D}.
\]
The map $q$ is continuous, and $M=q(L)$ is compact and
metrizable. If $q(y)=q(y')$, the continuous functions
$F(\cdot,y)$ and $F(\cdot,y')$ agree on the dense set $D$
and hence on all of $K$. Thus the function
\[
 \widetilde F:K\times M\to\R,\qquad
 \widetilde F(x,q(y))=F(x,y)
\]
is well defined. Its sections in the first variable are
continuous. For fixed $x$, continuity in the second variable
follows because the continuous surjection $q$ from a compact
space onto a Hausdorff space is a quotient map.
By Lemma~\ref{metric}, $\widetilde F\in B_1(K\times M)$.
Composition with $\operatorname{id}_K\times q$ completes the proof.
\end{proof}

\section{The counterexample}

\begin{lemma}\label{ccc}
The space $E$ satisfies the countable chain condition (ccc).
\end{lemma}
\begin{proof}
Let $\mathcal B$ be a Hamel basis of $\ell_\infty(\omega_1)$.
The map
\[
 J:E\to\R^{\mathcal B},\qquad J(a)=(L_b(a))_{b\in\mathcal B},
\]
is a topological embedding: the functionals $L_b$,
$b\in\mathcal B$, separate points and generate the given weak
topology, since every member of $\ell_\infty(\omega_1)$ is a
finite linear combination of members of $\mathcal B$.

The image $J(E)$ is dense in $\R^{\mathcal B}$.
Indeed, for distinct $b_1,\ldots,b_k\in\mathcal B$, the
functionals $L_{b_1},\ldots,L_{b_k}$ are linearly independent.
Consequently, the map
$a\mapsto(L_{b_1}(a),\ldots,L_{b_k}(a))$ is onto $\R^k$.
Otherwise a nonzero linear functional on $\R^k$ vanishing on
its image would give a nontrivial linear dependence among
the $L_{b_j}$.

For completeness, recall why a product of real lines is ccc.
Given an uncountable family of nonempty basic open cylinders,
the $\Delta$-system lemma yields an uncountable subfamily
whose finite supports have a common pairwise intersection,
called the root. Refining to rational intervals and passing
to an uncountable subfamily, we may assume that the
restrictions on the root are identical. Any two resulting
cylinders intersect. Thus $\R^{\mathcal B}$ is ccc.

Finally, a dense subspace of a ccc space is ccc. Indeed,
pairwise disjoint nonempty relatively open sets lift to
pairwise disjoint open sets of the ambient space: any
nonempty intersection of the latter would meet the dense
subspace. Apply this to $J(E)$.
\end{proof}

\begin{lemma}\label{cube}
The cube $K=[-1,1]^{\omega_1}$ is separable and nonmetrizable.
\end{lemma}
\begin{proof}
Choose distinct $t_\alpha\in(0,1)$.
For every polynomial $p$ with rational coefficients, put
\[
 x^p_\alpha=\max(-1,\min(1,p(t_\alpha))).
\]
These points form a countable dense subset of $K$.
Indeed, arbitrary prescribed values at finitely many distinct
labels can be interpolated by a real polynomial, whose
coefficients can then be approximated sufficiently closely
by rationals. Truncation to $[-1,1]$ does not increase the
error relative to a target value in $[-1,1]$.

The cube $K$ has no countable local base at zero.
Otherwise, choosing a basic cylinder about zero inside each
member of such a base would yield a countable set of
coordinates used by all these cylinders. A neighborhood
restricting a coordinate outside that set would contain none
of the chosen cylinders, a contradiction.
Thus $K$ is nonmetrizable.
\end{proof}

\begin{theorem}\label{main}
Assume $\bb>\omega_1$. The spaces
\[
 X=[-1,1]^{\omega_1},\qquad Y=E
\]
are completely regular and satisfy
\[
 \SC(X^2,\R)\subseteq B_1(X^2),\qquad
 \SC(Y^2,\R)\subseteq B_1(Y^2).
\]
However, the function
\[
 h:X\times Y\to\R,\qquad
 h(x,a)=\sum_{\alpha<\omega_1}x_\alpha a_\alpha
\]
is separately continuous and does not belong to $B_1(X\times Y)$.
\end{theorem}

\begin{proof}
Complete regularity follows from compactness of the Hausdorff
cube and from the Hausdorff locally convex topology of $E$.
The assertion concerning $Y^2$ is Corollary~\ref{Esquare}.
The assertion concerning $X^2$ follows from
Lemmas~\ref{compact} and~\ref{cube}.

For fixed $a$, the function $h(\cdot,a)$ is a finite linear
combination of coordinate projections. For fixed
$x\in X\subseteq\ell_\infty(\omega_1)$, the function
$h(x,\cdot)=L_x$ is continuous on $E$.
Hence $h$ is separately continuous.

Suppose that $h_n\in C(X\times E)$ and $h_n\to h$ pointwise.
For each $n$, consider
\[
 T_n:E\to C(X),\qquad T_n(a)=h_n(\cdot,a),
\]
where $C(X)$ carries the supremum norm.
The map $T_n$ is continuous. To see this, fix $a_0$ and
$\varepsilon>0$. Joint continuity provides, for each $x\in X$,
neighborhoods $U_x$ of $x$ and $V_x$ of $a_0$ such that
\[
 |h_n(z,a)-h_n(z,a_0)|<\varepsilon
 \quad(z\in U_x,\ a\in V_x).
\]
A finite subcover of $X$ by the $U_x$, together with the
intersection of the corresponding $V_x$, gives the desired
supremum norm estimate.

By Lemma~\ref{ccc}, $E$ is ccc, and so is its continuous
image $T_n(E)$.
Every metrizable ccc space is separable: for each positive
integer $k$, a maximal $1/k$-separated set is countable,
and the union of these sets is dense.
Choose a countable norm-dense set
$\mathcal D_n\subseteq T_n(E)$ and put
$\mathcal D=\bigcup_n\mathcal D_n$.

Suppose that $x,x'\in X$ satisfy $u(x)=u(x')$ for every
$u\in\mathcal D$. By uniform density, the equality holds
for all $u\in T_n(E)$ and all $n$. Consequently,
\[
 h_n(x,a)=h_n(x',a)\qquad(n\in\N,\ a\in E).
\]
Taking limits gives $h(x,a)=h(x',a)$ for every $a\in E$.
Substituting $a=e_\alpha$ yields $x_\alpha=x'_\alpha$.
Thus $x=x'$.

The countable family $\mathcal D\subseteq C(X)$ therefore
separates the points of $X$. The evaluation map
$X\to\R^{\mathcal D}$ is continuous and injective; compactness
makes it a homeomorphism onto its image.
This implies that $X$ is metrizable, contrary to
Lemma~\ref{cube}. Hence $h\notin B_1(X\times Y)$.
\end{proof}

\section*{AI best practices}
\addcontentsline{toc}{section}{AI best practices}
Codex assisted with the development and revision of the argument.
Parallel agents were used to examine intermediate proofs.

\end{document}